\documentclass[11pt]{article}

\newtheorem{thm}{Theorem}

\newtheorem{prop}{Proposition}
\newtheorem{conjecture}{Conjecture}
\newtheorem{defn}{Definition}[section]

\newcommand{\set}[1]{\left\{#1\right\}}
\newcommand{\td}{\mathrm{d}}
\newcommand{\error}{\mathcal{O}}

\begin{document}

{\LARGE On Partitions of Goldbach's Conjecture}

\bigskip

\noindent Max See Chin Woon

\noindent DPMMS, Cambridge University, Cambridge CB2 2BZ, UK

\smallskip

\noindent WatchMyPrice Technologies, 200 West Evelyn Ave, Suite
200, Mountain View, CA 94041, USA

\begin{quotation}\small
An approximate formula for the partitions of Goldbach's Conjecture
is derived using Prime Number Theorem and a probabilistic
approach.  A strong form of Goldbach's conjecture follows in the
form of a lower bounding function for the partitions of Goldbach's
conjecture.  Numerical computations suggest that the lower and
upper bounding functions for the partitions satisfy a simple
functional equation.  Assuming that this invariant scaling
property holds for all even integer $n$, the lower and upper
bounds can be expressed as simple exponentials.
\end{quotation}

\section{Goldbach's Conjecture and Recent Progress}

Goldbach's Conjecture states that every even integer $>2$ can be
expressed as a sum of two primes.

\bigskip

The proof remains an unsolved problem since Goldbach first wrote
the conjecture in a letter to Euler in 1792.  However, significant
progress has been made in recent years.

\bigskip

On the front of verifying Goldbach's Conjecture, no
counter-example has been found to date.  In 1855, Desboves
\cite{bib:Desboves1855} verified Goldbach's Conjecture for
$n<10,000$.  In 1940, Pipping \cite{bib:Pipping1938} verified the
conjecture for $n<100,000$ by the use of a computer. In 1964,
Stein \& Stein \cite{bib:Stein1965a, bib:Stein1965b} verified the
conjecture for $n<10^8$. In 1989, Granville, d. Lune, te Riele
\cite{bib:Granville1989} verified the conjecture for $n < 2 \times
10^{10}$. In 1998, Deshouillers, te Riele, Saouter
\cite{bib:Deshoulliers1998} verified the conjecture for $n <
10^{14}$, and Richstein \cite{bib:Richstein2000} verified the
conjecture for $n < 4 \times 10^{14}$.

\bigskip

On the front of proving {\em weaker forms} of Goldbach's
Conjecture, many have contributed.  In 1930, Shnirelman
\cite{bib:Shnirelman} achieved the first breakthrough by proving
that every natural number can be expressed as the sum of not more
than 20 primes.  In 1937, Vinagradov \cite{bib:Vinogradov1937}
proved that every sufficiently large odd numbers $>N_1$ can be
expressed as a sum of three primes.  In 1966, Chen
\cite{bib:Chen1966} proved that every sufficiently large even
natural number is the sum of a prime and a product of at most two
primes.  In 1995, Ramar\'e \cite{bib:Ramare1995, bib:Ramare1996}
proved that every even integer can be expressed as the sum of six
or fewer primes. In the same year, Kaniecki
\cite{bib:Kaniecki1995} proved that, assuming the Riemann
hypothesis is true, every odd integer can be expressed as the sum
of at most five primes. Kaniecki's proof can be improved to at
most four primes pending a further computational verification.

\bigskip

On the front of studying the properties of the partitions of
Goldbach's Conjecture, many asymptotic formulae have been derived
by Vinagradov \cite{bib:Vinogradov1937}, Chen \cite{bib:Chen1966,
bib:Chen1978} and others for partitions of the weaker forms of
Goldbach's Conjecture but none has been derived for partitions of
Goldbach's Conjecture.  In 1923, Hardy \& Littlewood
\cite{bib:HardyLittlewood1923} conjectured an asymptotic formula
for the partitions of Goldbach's Conjecture. The asymptotic
formula still remains to be proved to date.

\bigskip

In this paper, an approximate formula for the partitions of
Goldbach's Conjecture is derived using Prime Number Theorem and a
probabilistic approach.  In contrast to the many weaker forms of
Goldbach's Conjecture, a strong form of Goldbach's conjecture
follows in the form of a lower bounding function for the
partitions of Goldbach's conjecture.  The lower bounding function
can be estimated by the approximate formula which involves a sum
of the products of the reciprocals of logarithms.

\bigskip

Numerical computations suggest that the lower and upper bounding
functions for the partitions satisfy a simple functional equation.
Assuming that this invariant scaling property holds for all even
integer $n$, the lower and upper bounds can be expressed as simple
exponentials.

\section{Definitions}

\begin{defn}
Partition function for Goldbach's Conjecture, $G(n)$, is defined
as the number of representations of an even integer $n$ as the sum
of two primes $p$ and $q$.
\[ G(n) = \# \set{ (p,q) \mid n=p+q, p \le q }. \]
\end{defn}

\begin{prop}
Hardy and Littlewood \cite{bib:HardyLittlewood1923} gave a {\bf
conjectural} asymptotic formula,
\end{prop}

\begin{equation} \label{eqn:G(n)Asymptotic}
\lim_{n \to \infty} \frac{ G(n) }{ {\displaystyle \int_2^n
\frac{1}{(\ln x)^2} \, \td x \prod_{\stackrel{k=2}{p_k \mid n}}
\frac{p_k - 1}{p_k - 2} }} = 2 \; C
\end{equation}

where Hardy-Littlewood constant,

\begin{equation}
C = \prod_{k=2}^\infty \frac{ p_k (p_k - 2) }{ (p_k - 1)^2 } =
0.6601618158.
\end{equation}

\begin{thm}[Prime Number Theorem] \label{thm:PrimeNumberTheorem}
Let $\pi(n)$ be the number of primes $2 \le p \le n$.
\begin{equation}
\lim_{n \to \infty} \frac{\pi(n)}{n / \ln n} = 1.
\end{equation}
\end{thm}

\begin{thm}[Chebyshev's Limits] \label{thm:ChebyshevLimits}
Chebyshev \cite{bib:Rubinstein1994} established the limits such
that
\begin{equation}
\frac{7}{8} < \frac{\pi(n)}{n / \ln n} < \frac{9}{8}.
\end{equation}
\end{thm}

\section{Approximating $G(n)$ to first order}

The following theorem is derived based on Prime Number Theorem and
a heuristic probabilistic approach.

\begin{thm}
\begin{equation}
G(n) \approx \sum_{k=3}^{n/2} \frac{1}{\ln(k) \; \ln(n-k)}\;,
\quad \textrm{\rm even } n \ge 6. \label{eqn:G(n)Approx}
\end{equation}
\end{thm}

\noindent {\bf \it Proof.} From Theorem
\ref{thm:PrimeNumberTheorem} and \ref{thm:ChebyshevLimits}, the
Gauss's estimate $\pi(n) \approx n / \ln n$ is obtained.  The
probability of $n$ being a prime,
\begin{eqnarray}
P(n) & \approx & \frac{\pi(n) - \pi(n-1)}{n - (n-1)} \; = \; \pi(n) - \pi(n-1) \\
& \approx & \frac{\td}{\td x} \frac{x}{\ln x}\bigg\vert^{x=n} \; =
\; \frac{1}{\ln n} - \frac{1}{(\ln n)^2} \; \approx \;
 \frac{1}{\ln n} \label{eqn:P(n)}
\end{eqnarray}

\bigskip

Pick an integer $3 \le k \le n/2$, even $n \ge 6$.  The
probability of $k$ and $n-k$ both being primes is $P(k) \,
P(n-k)$.  Since the partition function $G(n)$ counts the number of
pairs of primes, $G(n)$ can be approximated by the sum of the
probabilities $P(k) \, P(n-k)$ over all k such that $3 \le k \le
n/2$.

\begin{eqnarray}
G(n) & \approx & \sum_{k = 3}^{n/2} P(k) \, P(n-k)\,, \quad {\rm
even } \; n \ge 6
\label{eqn:G(n)ProbabilityProductSum} \\
& \approx & \sum_{k=3}^{n/2} \frac{1}{\ln(k) \, \ln(n-k)}
\nonumber
\end{eqnarray}
\hfill Q.E.D.

\bigskip

While (\ref{eqn:G(n)Asymptotic}) is a more accurate approximation
than (\ref{eqn:G(n)Approx}), computation of
(\ref{eqn:G(n)Asymptotic}) requires the knowledge of the values of
all primes $p \mid n$. (\ref{eqn:G(n)Approx}) is a first order
approximation of $G(n)$ with no reference to the values of primes.

\section{A Strong Form of Goldbach's Conjecture}

From (\ref{eqn:P(n)}),
\begin{equation}
P(n) > \frac{1}{\ln n} - \left|\, \error\!\left(\frac{1}{(\ln
n)^2}\right)\right|
\end{equation}
It follows from (\ref{eqn:G(n)ProbabilityProductSum}) that
(\ref{eqn:G(n)Approx}) is approximating the lower bound of $G(n)$.
This leads to a conjecture stronger than Goldbach's Conjecture
where $G(n) > 0$ for even $n \ge 6$.

\begin{conjecture}[A Strong Form of Goldbach's Conjecture]
\begin{eqnarray}
G(n) & > & \sum_{k=3}^{n/2} \frac{1}{\ln(k) \;
\ln(n-k)} \nonumber \\
&& - \left|\, \error\!\!\left(\sum_{k=3}^{n/2}
\frac{1}{\left[\ln(k)\right] \; \left[\ln(n-k)\right]^2}
\right)\right| \nonumber \\
&& - \left|\, \error\!\!\left(\sum_{k=3}^{n/2}
\frac{1}{\left[\ln(k)\right]^2 \; \left[\ln(n-k)\right]}
\right)\right| \nonumber \\
&&  ( \textrm{\rm even } n \ge 6 ) \nonumber \\
& > & 0.
\end{eqnarray}
\end{conjecture}

\section{Lower and Upper Bounding Functions of $G(n)$}

\begin{defn}
The lower and upper bounding functions, respectively $f_L(x)$ and
$f_U(x)$, of $G(n)$ are defined as respectively monotonous
analytic functions such that $f_L(n) < G(n) < f_U(n)$ for even
integer $n \ge 6$.
\end{defn}

Numerical computations for $n<$ 1,000,000 suggest that $f_L(x)$
and $f_U(x)$ satisfy a functional equation of the form
\begin{equation}
\frac{f(ax)}{f(a)} = \frac{f(bx)}{f(b)}\;, \textrm{ constant }
a,b>0, \; x>0.
\end{equation}

A solution of the functional equation is
\begin{equation}
e^{\alpha \, x^\beta}, \textrm{ constant } \alpha>0, \; 0<\beta<1.
\end{equation}

This leads to the following conjecture.

\begin{conjecture}
The lower and upper bounding functions of $G(n)$ can be expressed
respectively as simple exponentials of the form $\exp(\alpha \,
x^{\,\beta})$, where constant $\alpha>0, \; 0<\beta<1$ can be
determined by numerical computations.
\end{conjecture}

\end{document}